%% file: Electromagnetics_Mirror.tex
\documentclass[]{interact}
\usepackage[T1]{fontenc}
\usepackage[utf8]{inputenc}
\usepackage{lmodern}
\usepackage{bm,mathtools,mathrsfs}
\usepackage{tabularx,array}
\usepackage{placeins}
\usepackage{needspace}
\usepackage{natbib}
\bibpunct[, ]{(}{)}{;}{a}{}{,}
\renewcommand\bibfont{\fontsize{10}{12}\selectfont}
\usepackage[noperiod]{jabbrv}
\usepackage{url}
\theoremstyle{plain}
\newtheorem{theorem}{Theorem}[section]

\theoremstyle{remark}

\theoremstyle{definition}

\newcommand{\RR}{\mathbb R}
\newcommand{\CC}{\mathbb C}
\newcommand{\EE}{\mathbf E}
\newcommand{\HH}{\mathbf H}
\newcommand{\JJ}{\mathbf J}
\newcommand{\nn}{\mathbf n}
\newcommand{\vv}{\mathbf v}
\newcommand{\xx}{\mathbf x}
\newcommand{\curl}{\operatorname{curl}}
\newcommand{\diver}{\operatorname{div}}
\newcommand{\I}{\mathbf I}
\newcommand{\TT}{\mathcal T}
\newcommand{\OO}{\mathrm O}
\newcommand{\eps}{\boldsymbol\varepsilon}
\newcommand{\mum}{\boldsymbol\mu}
\newcommand{\sig}{\boldsymbol\sigma}
\newcommand{\cL}{\mathcal L}
\newcommand{\Med}{\mathcal M}
\newcommand{\Bt}{\mathcal B_{\rm PEC}}

\begin{document}
\title{Symmetry Inheritance and Exact Mirror-Symmetry Reduction for Electric-Field Vector Wave Problems with PEC Boundaries}
\author{
\name{Wei Jiang\thanks{CONTACT Wei Jiang. Email: jwmathphy@163.com}}
\affil{School of Physics and Mechatronics Engineering, Guizhou Minzu University, Guiyang 550025, China}
}
\maketitle
\begin{abstract}
A symmetry-inheritance and exact domain-reduction framework is developed
for the second-order vector wave equation of the electric field in bounded
domains with a perfect electric conductor (PEC) boundary and general
spatially varying anisotropic media. Orthogonal group actions on polar
vectors, axial vectors, and second-order constitutive tensors are used to
derive the transformed curl--curl equation and PEC trace directly from the
original boundary-value problem. Under unique solvability, the common
symmetry group of the problem data is inherited by the electric field.
When the symmetry group of the electromagnetic medium admits a mirror symmetry, an arbitrary
impressed current can be decomposed by the associated reflection
projectors into mirror-even and mirror-odd components. The corresponding
full-domain fields inherit the respective parities, and their exact
restriction to one half-domain yields a perfect magnetic conductor (PMC)
condition for the even field and a PEC condition for the odd field. The
original problem is therefore represented exactly by two independent
reduced-domain problems. A symmetry-compatible first-family N\'ed\'elec
discretization preserves this decomposition at the discrete level.
Analytical and numerical results verify the symmetry-inheritance theory,
the exact PEC/PMC reduced-domain formulations, and the discrete
reconstruction for constant and spatially varying anisotropic conductive
media under both symmetric and asymmetric excitations.
\end{abstract}

\section{Introduction}\label{sec:intro}
The second-order electric-field vector wave equation is one of the standard formulations for time-harmonic electromagnetic analysis. It is particularly attractive for finite-element computation because the electric field can be approximated directly by curl-conforming edge elements and the perfect electric conductor condition is imposed through the tangential electric trace. The mathematical and computational foundations of this formulation are well established in the literature on edge elements \citep{monk2003,hiptmair2002}. Detailed treatments of practical electromagnetic finite-element modeling and edge-element formulations can also be found in \citep{jin2014,bossavit1998}.

Symmetry has long been used to reduce electromagnetic computations, but a geometrically symmetric domain does not by itself guarantee a symmetric solution. The governing operator, material tensors, physical boundary conditions, and excitation must transform consistently. Bossavit developed a systematic theory of boundary-value problems with geometrical symmetry and their finite-element approximation \citep{bossavit1986,bossavit1993}. Finite-group representation theory was introduced into field computation by Ballisti, Hafner, and Leuchtmann, and symmetry-adapted decompositions were later developed directly for Maxwell problems on symmetric structures \citep{ballisti1982,weiland2002}. Symmetry planes have also been used to reduce the number of unknowns in integral-equation computations of electromagnetic scattering and radiation \citep{johnson1990}. A related symmetry-inheritance strategy has recently been formulated for second-order linear elliptic problems, where uniqueness transfers the common symmetry of the problem data to the solution \citep{jiangpan2026robin}. The present work addresses the corresponding vector electric-field problem and the additional complications created by curl operators and tangential traces under reflections.

The role of uniqueness in the present analysis is deliberately separated from the symmetry argument itself. Classical Maxwell theory provides the functional setting and well-posedness results for bounded electromagnetic boundary-value problems \citep{cessenat1996,costabel1991}. Compactness and Fredholm properties, together with exclusion of resonant cavity frequencies in the lossless case, provide the analytical basis for the uniqueness assumption used below \citep{weber1980,monkdemkowicz2001}. The symmetry theorem therefore does not require a new existence theory and it requires only unique solvability of the corresponding perfect electric conductor problem at the operating frequency.

A correct reduced formulation also depends on the boundary-trace structure of curl-conforming fields. Tangential trace spaces and the integration-by-parts identities needed for Maxwell problems on Lipschitz domains and polyhedra have been studied in detail \citep{buffaciarlet2001a,buffacostabel2002}. The classical first-family N\'ed\'elec elements provide the conforming discretization adopted in this paper \citep{nedelec1980,nedelec1986}. Compatible and differential-complex viewpoints explain why edge orientations and tangential degrees of freedom are essential in electromagnetic discretizations \citep{arnold2006,arnold2010}. Closely related structure-preserving approaches include the finite integration technique and geometric simplicial discretizations of Maxwell equations \citep{clemens2001,teixeira2001}.

The treatment of material invariance under spatial transformations is
especially important for anisotropic media. Under a spatial reflection,
the constitutive parameters transform according to the transformation law
for second-order tensor fields. Therefore, invariance of an anisotropic
medium under reflection is determined by the transformed tensor field as a
whole and cannot, in general, be reduced to entrywise symmetry of the
constitutive matrices. This notion of spatial invariance is distinct from
algebraic properties of the constitutive tensors, such as matrix symmetry
or Hermitian symmetry. Standard electromagnetic treatments of
inhomogeneous and anisotropic media emphasize the tensorial nature of the
constitutive relations \citep{chew1995,kong1990}. The influence of
anisotropy on electromagnetic modal structure has also been studied in
waveguide problems, where relations between the constitutive tensors
determine whether pure transverse-electric and transverse-magnetic modes
can exist independently \citep{jiang2017,jiang2018}. These considerations
motivate treating the spatial transformations of the permittivity,
permeability, and conductivity tensors explicitly throughout the present
analysis.

The central problem considered here is source-induced breaking of a mirror symmetry. The geometry, exterior perfect electric conductor boundary, and constitutive medium may all possess a reflection symmetry while the impressed current has no definite mirror parity. In that situation the total field cannot in general be reduced by imposing a perfect electric conductor or perfect magnetic conductor condition directly on the symmetry plane. The present work first decomposes the arbitrary source into its mirror-even and mirror-odd components on the original full domain. Unique solvability then forces the two associated electric fields to inherit the corresponding parities. Only after this full-domain result has been established are the two problems restricted to one half-domain, where exact perfect magnetic conductor and perfect electric conductor internal conditions follow from Maxwell's equations.

The principal contributions are therefore threefold. First, the transformation of the original second-order electric-field boundary-value problem is derived by distinguishing polar-vector, axial-vector, and tensor actions, including the orientation change of curl and the perfect electric conductor trace under reflections. Second, an exact mirror decomposition is established for arbitrary impressed currents, including the practically important case in which the excitation destroys the mirror symmetry of the total field. Third, the resulting half-domain problems are transferred to a N\'ed\'elec discretization with globally oriented edge degrees of freedom and are verified on analytical and three-dimensional numerical examples, including an anisotropically loaded shielded cavity with a localized asymmetric excitation. The detailed chain-rule, curl, tensor, normal-vector, and boundary-trace calculations are collected in the Appendix so that the main text can focus on the symmetry structure, exact internal boundary conditions, and computational consequences.

\section{Electric-field vector boundary-value problem and uniqueness}\label{sec:model}
Let $\Omega\subset\RR^3$ be a bounded, connected Lipschitz domain with outward unit normal $\nn$, and assume that the entire boundary $\partial\Omega$ is a PEC. The time-harmonic convention $e^{j\omega t}$ is adopted, with $\omega>0$. The constitutive parameters
$\eps,\mum,\sig:\Omega\rightarrow\CC^{3\times3}$
are allowed to be general spatially varying tensor fields. The permeability tensor $\mum(\xx)$ is assumed nonsingular in $\Omega$ so that $\mum(\xx)^{-1}$ is well defined. No Hermitian, positivity, reciprocity, isotropy, or losslessness assumptions are required for the symmetry analysis. Sufficient regularity is assumed whenever pointwise differential identities are used. When symmetry inheritance is invoked, the corresponding PEC boundary-value problem is assumed to admit a unique
solution.

Collect the constitutive tensors into the medium data
\[
\Med:=(\eps,\mum,\sig),
\]
and define the corresponding vector wave operator by
\begin{equation}\label{eq:Aomega}
\cL_{\Med,\omega}\mathbf v
:=\curl\!\left(\mum^{-1}\curl\mathbf v\right)
-\omega^2\eps\mathbf v+j\omega\sig\mathbf v.
\end{equation}
The boundary-value problem is
\begin{subequations}\label{eq:model}
\begin{align}
\cL_{\Med,\omega}\EE
&=-j\omega\JJ &&\text{in }\Omega,\label{eq:modela}\\
\Bt\EE:=\nn\times\EE
&=\mathbf0 &&\text{on }\partial\Omega.\label{eq:modelb}
\end{align}
\end{subequations}
Here $\eps$ and $\sig$ are the permittivity and electrical conductivity tensors, respectively, and $\JJ$ is the impressed electric-current density. Both $\EE$ and $\JJ$ are polar vector fields. The form of \eqref{eq:modela} follows directly from the frequency-domain Maxwell equations with the stated time convention.

For the symmetry analysis, the essential assumption is uniqueness. Under the standard Maxwell Fredholm setting  \eqref{eq:model} is uniquely solvable precisely when the corresponding homogeneous problem
\begin{subequations}\label{eq:hom}
\begin{align}
\cL_{\Med,\omega}\mathbf u&=\mathbf0 &&\text{in }\Omega,\\
\nn\times\mathbf u&=\mathbf0 &&\text{on }\partial\Omega
\end{align}
\end{subequations}
has only the zero solution. A nonzero homogeneous field destroys uniqueness immediately, while a trivial kernel implies invertibility when the operator is Fredholm of index zero \citep{cessenat1996,weber1980}.

In the lossless case $\sig=\mathbf0$, nonzero solutions of \eqref{eq:hom} are PEC cavity modes. On the usual divergence-compatible subspace they satisfy
\begin{subequations}\label{eq:eig}
\begin{align}
\curl(\mum^{-1}\curl\mathbf u)&=\lambda\eps\mathbf u &&\text{in }\Omega,\\
\diver(\eps\mathbf u)&=0&&\text{in }\Omega,\\
\nn\times\mathbf u&=\mathbf0 &&\text{on }\partial\Omega,
\end{align}
\end{subequations}
with $\lambda=\omega^2$. Thus a lossless driven problem is unique whenever the working frequency is not a PEC Maxwell eigenfrequency. In conductive media the same symmetry theory applies whenever \eqref{eq:hom} has a trivial kernel. For example, a uniformly positive-definite conductivity tensor is a standard sufficient condition at a real nonzero frequency under additional Hermitian and positivity assumptions. 

\section{Symmetry groups and symmetry inheritance}\label{sec:symmetry}
Choose the origin at a common fixed point of the spatial symmetries and define the isometry group of the bounded domain $\Omega$ by
\[
G_\Omega:=\{g\in\OO(3):g\Omega=\Omega\},
\]
where $\OO(3)$ is the full orthogonal group in three dimensions. Translations are excluded here because $\Omega$ is bounded, and a nonzero translation cannot map a bounded domain onto itself.

For $g\in G_\Omega$, let $\TT_g^0$, $\TT_g^p$, $\TT_g^a$, and $\TT_g^2$ denote the left group actions on a scalar field $f$, a polar vector $\mathbf v$, an axial vector $\mathbf w$, and a second-order tensor $\mathbf A$, respectively:
\begin{equation}\label{eq:mainactions}
\begin{aligned}
(\TT_g^0 f)(\xx)&= f(g^{-1}\xx),\\
(\TT_g^p\mathbf v)(\xx)&=g\mathbf v(g^{-1}\xx),\\
(\TT_g^a\mathbf w)(\xx)&=\det(g)\,g\mathbf w(g^{-1}\xx),\\
(\TT_g^2\mathbf A)(\xx)&=g\mathbf A(g^{-1}\xx)g^{-1}.
\end{aligned}
\end{equation}
The distinction between $\TT_g^p$ and $\TT_g^a$ is essential for reflections. In particular, Appendix~\ref{app:actions} proves
\begin{equation}\label{eq:mainintertwine}
\curl(\TT_g^p\mathbf v)=\TT_g^a(\curl\mathbf v),
\qquad
\curl(\TT_g^a\mathbf w)=\TT_g^p(\curl\mathbf w),
\end{equation}
and the PEC-trace relation
\begin{equation}\label{eq:pecmain}
\Bt(\TT_g^p\mathbf v)=\TT_g^a(\Bt\mathbf v),
\qquad g\in G_\Omega.
\end{equation}
The complete derivations, including the transformation of the outward normal and the curl--curl term, are given in Appendix~\ref{app:actions}.

Set
\begin{gather}
\eps_g=\TT_g^2\eps, \qquad
\mum_g=\TT_g^2\mum, \qquad
\sig_g=\TT_g^2, \qquad
\JJ_g=\TT_g^p\JJ,\qquad
\EE_g=\TT_g^p\EE,
\end{gather}
and denote the transformed medium by
\[
\Med_g:=(\eps_g,\mum_g,\sig_g).
\]
The covariance identity proved in Appendix~\ref{app:actions} gives the transformed problem
\begin{subequations}\label{eq:transformed}
\begin{align}
\cL_{\Med_g,\omega}\EE_g
&=-j\omega\JJ_g &&\text{in }\Omega,\label{eq:transformed-a}\\
\Bt\EE_g&=\mathbf0 &&\text{on }\partial\Omega.\label{eq:transformed-b}
\end{align}
\end{subequations}
Thus the strong partial differential equation \eqref{eq:transformed-a} and its homogeneous PEC trace \eqref{eq:transformed-b} transform together.

The constitutive tensor and source symmetry groups are
\begin{align*}
G_\varepsilon&:=\{g\in G_\Omega:\TT_g^2\eps=\eps\},\\
G_\mu&:=\{g\in G_\Omega:\TT_g^2\mum=\mum\},\\
G_\sigma&:=\{g\in G_\Omega:\TT_g^2\sig=\sig\},\\
G_J&:=\{g\in G_\Omega:\TT_g^p\JJ=\JJ\}.
\end{align*}
Define the symmetry group of the electromagnetic medium 
$$G_{\rm med}:=G_\varepsilon\cap G_\mu\cap G_\sigma.$$ Since \eqref{eq:pecmain} shows that the homogeneous PEC condition on the whole physical boundary is preserved by every $g\in G_\Omega$, its symmetry group is simply $G_{\rm PEC}=G_\Omega$. Consequently, the common symmetry group of all data defining \eqref{eq:model} is
\begin{equation}\label{eq:Gdata}
G_{\rm data}:=G_\varepsilon\cap G_\mu\cap G_\sigma\cap G_J.
\end{equation}
The symmetry group of the solution to \eqref{eq:model} is
\[
G_E:=\{g\in G_\Omega:\TT_g^p\EE=\EE\}.
\]

\begin{theorem}\label{thm:inherit}
Assume that \eqref{eq:model} is uniquely solvable. Then
\begin{equation}\label{eq:subgroup}
G_{\rm data}\subseteq G_E.
\end{equation}
\end{theorem}
\begin{proof}
Let $g\in G_{\rm data}$ and set $\EE_g=\TT_g^p\EE$. By definition,
$\eps_g=\eps$, $\mum_g=\mum$, $\sig_g=\sig$, and $\JJ_g=\JJ$. Hence \eqref{eq:transformed-a} is exactly \eqref{eq:modela}. Equation~\eqref{eq:pecmain} also gives $\Bt\EE_g=0$ on $\partial\Omega$. Therefore $\EE_g$ and $\EE$ solve the same complete PEC boundary-value problem. Uniqueness implies $\EE_g=\EE$, so $g\in G_E$. Since $g$ is arbitrary, \eqref{eq:subgroup} follows.
\end{proof}

A particularly simple example illustrates why the reverse inclusion
$G_E\subseteq G_{\rm data}$ cannot be asserted in general. Consider a
spherical PEC cavity whose material distribution has no nontrivial
orthogonal symmetry, so that $G_{\rm data}=\{I\}$. Let the impressed current be identically zero. The problem then becomes the corresponding homogeneous PEC problem.
Assume that the operating frequency is such that the homogeneous problem
has no nonzero solution. By uniqueness, $\EE=\mathbf0$ holds. Since the zero field is invariant under every orthogonal transformation
and the sphere has the full symmetry group $O(3)$, the solution satisfies $G_E=O(3)$. 
Consequently,
\[
G_{\rm data}=\{I\}\subsetneq O(3)=G_E.
\]

This simple example shows that a particular solution may possess
additional symmetries that are not shared by the partial differential equation (PDE) data. Hence, only
the inclusion
$G_{\rm data}\subseteq G_E$
can be guaranteed in general. $G_{\rm data}$ represents the symmetries guaranteed from the PDE data before the field is computed.

\section{Mirror decomposition and exact internal boundary conditions}\label{sec:mirror}
We take the volume centroid of the bounded domain $\Omega$ as the
coordinate origin. Since every Euclidean isometry that maps $\Omega$
onto itself fixes its centroid, the isometry group of $\Omega$ is
identified, in these coordinates, with a subgroup of $\OO(3)$.
Now consider reflection in a plane through this origin. Let
\[
\Pi:=\{\xx\in\RR^3:\nn_S\cdot\xx=0\},
\qquad
R:=\I-2\nn_S\nn_S^T,
\]
where $\nn_S$ is a unit normal on $\Pi$. Then $R$ is the reflection in $\Pi$, with $R^T=R$, $R^2=\I$ and $\det R=-1$. Assume that $G_{\rm med}$ contains this reflection, i.e.,
\begin{equation}\label{eq:bgmirror}
R\in G_\varepsilon\cap G_\mu\cap G_\sigma.
\end{equation}
No mirror-parity assumption is imposed on the impressed current $\JJ$. In particular, the generic source-symmetry-breaking case considered below satisfies
\[
\TT_R^p\JJ\neq\JJ,
\qquad
\TT_R^p\JJ\neq-\JJ,
\]
so that the total field has no definite mirror parity that can be inferred directly from Theorem~\ref{thm:inherit}. The decomposition that follows, however, is valid for an arbitrary impressed current.

Write $T:=\TT_R^p$. Since $T^2=I$, the two reflection projectors
\begin{equation}\label{eq:proj}
P_+:=\frac12(I+T),
\qquad
P_-:=\frac12(I-T)
\end{equation}
are complementary. Decompose the impressed current $\JJ$ as
\begin{equation}\label{eq:Jpm}
\JJ=\JJ_++\JJ_-,
\end{equation}
where 
\begin{equation*}
\JJ_{+}=P_{+}\JJ,\qquad
\JJ_{-}=P_{-}\JJ.
\end{equation*}
Clearly, 
\begin{equation*}
T\JJ_{+}=\JJ_{+},\qquad
T\JJ_{-}=-\JJ_{-}.
\end{equation*}
The two associated electric fields are now defined on the entire original domain. Let $\EE_+$ and $\EE_-$ be the unique solutions of
\begin{subequations}\label{eq:fullplus}
\begin{align}
\cL_{\Med,\omega}\EE_+
&=-j\omega\JJ_+ &&\text{in }\Omega,\label{eq:fullplus-a}\\
\Bt\EE_+&=\mathbf0 &&\text{on }\partial\Omega,
\label{eq:fullplus-b}
\end{align}
\end{subequations}
and
\begin{subequations}\label{eq:fullminus}
\begin{align}
\cL_{\Med,\omega}\EE_-
&=-j\omega\JJ_- &&\text{in }\Omega,\label{eq:fullminus-a}\\
\Bt\EE_-&=\mathbf0 &&\text{on }\partial\Omega.
\label{eq:fullminus-b}
\end{align}
\end{subequations}
Because the differential operator, domain, and PEC exterior boundary are invariant under $R$, applying $T$ to \eqref{eq:fullplus} produces another solution driven by $\JJ_+$. Uniqueness gives
\[
T\EE_+=\EE_+.
\]
For \eqref{eq:fullminus}, the transformed field $T\EE_-$ is driven by $-\JJ_-$. The field $-\EE_-$ satisfies exactly the same transformed problem, therefore uniqueness gives
\[
T\EE_-=-\EE_-.
\]
Thus the two full-domain fields have definite electric parities. Moreover, linearity shows that $\EE_++\EE_-$ solves \eqref{eq:model} with source $\JJ_++\JJ_-=\JJ$, and uniqueness yields
\begin{equation}\label{eq:superposition}
\EE=\EE_++\EE_-.
\end{equation}
This full-domain decomposition is the starting point for the subsequent region reduction.

We next derive the conditions induced on the fixed plane $\Pi$. At $\xx\in\Pi$, $R\xx=\xx$, and every vector has the orthogonal decomposition
\[
\mathbf v=\mathbf v_t+v_n\nn_S,
\qquad
R\mathbf v=\mathbf v_t-v_n\nn_S.
\]

For the odd field, $T\EE_-=-\EE_-$ gives $R\EE_{-}(\xx)=-\EE_{-}(\xx)$ for $\xx\in{\Pi}$. Hence $(\EE_-)_{t}=0$ and
\begin{equation}\label{eq:oddpec}
\nn_S\times\EE_-=\mathbf0
\qquad\text{on }\Pi\cap\Omega.
\end{equation}
This is the electric condition of a PEC symmetry plane. The axial-vector curl transformation also gives $R\curl\EE_-=\curl\EE_-$ on $\Pi$, so $\curl\EE_-$ is tangential and Faraday's law implies $\nn_S\cdot\mathbf B_-=0$. Thus the usual PEC pair is recovered.

For the even field, $T\EE_+=\EE_+$ implies, at every fixed point $\xx\in\Pi$,
\begin{equation}\label{eert}
R\EE_+(\xx)=\EE_+(\xx).
\end{equation}
Hence the normal component of $\EE_+$ vanishes and $\EE_+$ is tangential to the mirror plane. Applying the first identity in \eqref{eq:mainintertwine} gives
\[
\curl\EE_+=\curl(T\EE_+)=\TT_R^a(\curl\EE_+).
\]
Because $\det R=-1$, the axial action satisfies $(\TT_R^a\mathbf w)(\xx)=-R\mathbf w(\xx)$ at a fixed point of the reflection. Therefore,
\begin{equation}\label{curlre}
R\curl\EE_+(\xx)=-\curl\EE_+(\xx),
\qquad \xx\in\Pi,
\end{equation}
so $\curl\EE_+$ is normal to $\Pi$. The permeability symmetry in \eqref{eq:bgmirror} gives
\[
R\mum(\xx)R^T=\mum(\xx),
\qquad \xx\in\Pi,
\]
and, since $R^T=R$, this is equivalent to $R\mum^{-1}=\mum^{-1}R$ on the mirror plane. Combining this commutation relation with \eqref{curlre} yields
\begin{equation}\label{curlret}
R\mum^{-1}\curl\EE_+
=-\mum^{-1}\curl\EE_+
\qquad\text{on }\Pi.
\end{equation}
Thus $\mum^{-1}\curl\EE_+$ is normal to the mirror plane and
\begin{equation}\label{eq:evenpmc}
\nn_S\times\mum^{-1}\curl\EE_+=\mathbf0
\qquad\text{on }\Pi\cap\Omega.
\end{equation}
With the adopted time convention, $\mum^{-1}\curl\EE=-j\omega\HH$, so \eqref{eq:evenpmc} is equivalent to $\nn_S\times\HH_+=0$. Likewise, the permittivity symmetry gives $R\eps=\eps R$ on $\Pi$. Together with \eqref{eert}, this implies
\[
R(\eps\EE_+)=\eps\EE_+
\qquad\text{on }\Pi,
\]
so $\eps\EE_+$ is tangential to the mirror plane and therefore $\nn_S\cdot\mathbf D_+=0$ for $\mathbf D_+=\eps\EE_+$. Thus the usual PMC pair is recovered.

Let $\Omega_+$ be one half of $\Omega$ and write
\[
\Gamma_P:=\partial\Omega_+\cap\partial\Omega,
\qquad
\Gamma_S:=\partial\Omega_+\cap\Pi.
\]
Restricting the two full-domain fields to $\Omega_+$ yields the exact reduced problems
\begin{subequations}\label{eq:plusstrong}
\begin{align}
\cL_{\Med,\omega}\mathbf u_+
&=-j\omega\JJ_+ &&\text{in }\Omega_+,\\
\nn\times\mathbf u_+&=\mathbf0 &&\text{on }\Gamma_P,\\
\nn_S\times\mum^{-1}\curl\mathbf u_+&=\mathbf0 &&\text{on }\Gamma_S,
\end{align}
\end{subequations}
and
\begin{subequations}\label{eq:minusstrong}
\begin{align}
\cL_{\Med,\omega}\mathbf u_-
&=-j\omega\JJ_- &&\text{in }\Omega_+,\\
\nn\times\mathbf u_-&=\mathbf0 &&\text{on }\Gamma_P,\\
\nn_S\times\mathbf u_-&=\mathbf0 &&\text{on }\Gamma_S.
\end{align}
\end{subequations}
Thus the even problem \eqref{eq:plusstrong} has a PMC artificial boundary and the odd problem \eqref{eq:minusstrong} has a PEC artificial boundary. Their parity extensions are
\[
(\mathscr E_+\mathbf u_+)(\xx)=
\begin{cases}
\mathbf u_+(\xx),&\xx\in\Omega_+,\\
R\mathbf u_+(R\xx),&\xx\in\Omega_-,
\end{cases}
\]
and
\[
(\mathscr E_-\mathbf u_-)(\xx)=
\begin{cases}
\mathbf u_-(\xx),&\xx\in\Omega_+,\\
-R\mathbf u_-(R\xx),&\xx\in\Omega_-,
\end{cases}
\]
where $\Omega_-=\Omega\backslash\overline{\Omega_{+}}$.
Under the stated regularity, the interface conditions above are exactly those required for the reflected fields to satisfy the corresponding full-domain problems. Consequently,
\begin{equation}\label{eq:reconstruct}
\EE=\mathscr E_+\mathbf u_++\mathscr E_-\mathbf u_-.
\end{equation}
For a generic asymmetric source both half-domain problems must be solved. If the source already has definite mirror parity, one projector vanishes and only one half-domain solve is required.

\begin{theorem}\label{thm:mirrorreduction}
Assume that the full-domain PEC problem \eqref{eq:model} is uniquely solvable and that the reflection $R$ belongs to $G_{\rm med}$. For an arbitrary impressed current $\JJ$, define $\JJ_+$ and $\JJ_-$ by \eqref{eq:proj}--\eqref{eq:Jpm}. Then the unique full-domain solution admits the exact decomposition
\[
\EE=\EE_++\EE_-,
\qquad
\TT_R^p\EE_+=\EE_+,
\qquad
\TT_R^p\EE_-=-\EE_-.
\]
Moreover, the restrictions $\mathbf u_+=\EE_+|_{\Omega_+}$ and $\mathbf u_-=\EE_-|_{\Omega_+}$ are the unique solutions of the reduced problems \eqref{eq:plusstrong} and \eqref{eq:minusstrong}, respectively, and the original field is recovered exactly by \eqref{eq:reconstruct}.
\end{theorem}
\begin{proof}
The full-domain parity relations and the superposition identity follow from the reflection projectors, invariance of the background problem, linearity, and uniqueness, as shown above. Equations~\eqref{eq:oddpec} and \eqref{eq:evenpmc} give the exact internal PEC and PMC conditions, so the restrictions of $\EE_-$ and $\EE_+$ satisfy \eqref{eq:minusstrong} and \eqref{eq:plusstrong}. Conversely, the corresponding parity extension of any reduced homogeneous solution satisfies the full-domain homogeneous PEC problem. Full-domain uniqueness therefore implies uniqueness of both reduced problems. Applying the two parity extensions to their unique solutions and using \eqref{eq:superposition} gives \eqref{eq:reconstruct}.
\end{proof}

\section{N\'ed\'elec edge-element discretization}
\label{sec:discretization}

The reduced problems \eqref{eq:plusstrong} and \eqref{eq:minusstrong}
are discretized only after the exact PEC and PMC conditions on the
symmetry plane have been established at the continuous level. Define
\begin{align*}
V_+
&:=
\bigl\{
\vv\in H(\curl,\Omega_+):
\nn\times\vv=\mathbf0
\ \text{on }\Gamma_P
\bigr\},\\
V_-
&:=
\bigl\{
\vv\in H(\curl,\Omega_+):
\nn\times\vv=\mathbf0
\ \text{on }\partial\Omega_{+}
\bigr\}.
\end{align*}
Thus, the physical PEC condition on $\Gamma_P$ is imposed essentially
in both spaces. On the artificial boundary $\Gamma_S$, the odd problem
has the essential PEC condition, whereas the PMC condition of the even
problem is natural.

For $\mathbf u,\vv\in H(\curl,\Omega_+)$, define the sesquilinear form
\begin{align*}
a_\omega(\mathbf u,\vv)
:={}&
\int_{\Omega_+}
\bigl(\mum^{-1}\curl\mathbf u\bigr)
\cdot\overline{\curl\vv}\,d\xx
-\omega^2
\int_{\Omega_+}
(\eps\mathbf u)\cdot\overline{\vv}\,d\xx\\
&\quad
+j\omega
\int_{\Omega_+}
(\sig\mathbf u)\cdot\overline{\vv}\,d\xx ,
\end{align*}
and the linear functionals
\[
\ell_\pm(\vv)
:=
-j\omega
\int_{\Omega_+}
\JJ_\pm\cdot\overline{\vv}\,d\xx .
\]
For the even problem, integration by parts of the curl--curl term
produces the boundary trace
$\nn_S\times\mum^{-1}\curl\mathbf u_+$ on $\Gamma_S$, which vanishes
by \eqref{eq:evenpmc}. For the odd problem, the tangential traces of
the trial and test functions vanish on $\Gamma_S$ by
\eqref{eq:oddpec}. Consequently, the two reduced problems can be written
in the unified variational form
\begin{equation}
\begin{aligned}
&\text{find }\mathbf u_\pm\in V_\pm\text{ such that}\\
&a_\omega(\mathbf u_\pm,\vv)
=\ell_\pm(\vv)
\qquad
\forall\,\vv\in V_\pm .
\end{aligned}
\label{eq:weakpm}
\end{equation}

Let $V_{h,\pm}\subset V_\pm$ be the lowest-order first-family
N\'ed\'elec edge-element spaces \citep{nedelec1980,nedelec1986}, and let
$\{\mathbf N_j^\pm\}_{j=1}^{N_\pm}$ denote their global basis functions.
The discrete fields are written as
\[
\mathbf u_{h,\pm}
=
\sum_{j=1}^{N_\pm}
U_j^\pm\,\mathbf N_j^\pm .
\]
Substitution into \eqref{eq:weakpm} with
$\vv=\mathbf N_i^\pm$ gives the algebraic systems
\begin{equation}
\mathbf A_\pm\mathbf U_\pm
=
\mathbf b_\pm ,
\label{eq:algebraicpm}
\end{equation}
where
\[
\mathbf A_\pm
=
\mathbf K_\pm
-\omega^2\mathbf M_{\varepsilon,\pm}
+j\omega\mathbf M_{\sigma,\pm},
\]
and
\begin{align*}
(\mathbf K_\pm)_{ij}
&=
\int_{\Omega_+}
\bigl(\mum^{-1}\curl\mathbf N_j^\pm\bigr)
\cdot
\overline{\curl\mathbf N_i^\pm}\,d\xx,\\
(\mathbf M_{\varepsilon,\pm})_{ij}
&=
\int_{\Omega_+}
(\eps\mathbf N_j^\pm)
\cdot
\overline{\mathbf N_i^\pm}\,d\xx,\\
(\mathbf M_{\sigma,\pm})_{ij}
&=
\int_{\Omega_+}
(\sig\mathbf N_j^\pm)
\cdot
\overline{\mathbf N_i^\pm}\,d\xx,\\
(\mathbf b_\pm)_i
&=
-j\omega
\int_{\Omega_+}
\JJ_\pm\cdot
\overline{\mathbf N_i^\pm}\,d\xx .
\end{align*}
The standard edge degree of freedom is
\[
d_e(\vv)=\int_e\vv\cdot\mathbf t_e\,ds,
\]
with a fixed global orientation $\mathbf t_e$ assigned to every mesh
edge. Accordingly, the edge unknowns on $\Gamma_S$ are eliminated in
the odd problem, while they remain active in the even problem because
the PMC condition is natural. The full-domain fields are subsequently
reconstructed from the two half-domain solutions by the parity
extensions $\mathscr E_+$ and $\mathscr E_-$ defined in
\eqref{eq:reconstruct}, with the reflected edge orientations accounted
for by the prescribed global edge orientation.

\section{Analytical and numerical verification}\label{sec:numerics}
The analytical example below first verifies the symmetry-inheritance result independently of any discretization. The subsequent three examples are designed to test three different aspects of the proposed reduction: the exact PEC/PMC symmetry boundaries for fields of definite parity, tensorial mirror symmetry in a spatially varying anisotropic conductive medium, and the computational reduction of an engineering-motivated cavity problem driven by a genuinely asymmetric localized source.

For Examples~2 and~3, reflection-paired tetrahedral meshes of $\Omega=(-1,1)^3$ are used. The mesh on $z\geq0$ is generated first and then reflected across $z=0$, so that the two halves match exactly. Example~4 uses the same paired-mesh construction for a rectangular shielded cavity. Every global edge is oriented from the smaller to the larger global node number. The element matrices and load vectors are evaluated by an eleven-point degree-four tetrahedral quadrature rule, and all complex sparse systems are solved by sparse LU factorization. Let $N_{\rm full}$, $N_+$, and $N_-$ denote the active edge unknowns of the full-domain, electric-even, and electric-odd problems, respectively.

For the comparison between a direct full-domain solution and the field reconstructed from the two reduced problems, define
\begin{equation}\label{eq:deltarecnum}
\delta_{\rm rec}:=
\frac{\|\EE_h^{\rm full}-(\mathscr E_+\mathbf u_{h,+}+\mathscr E_-\mathbf u_{h,-})\|_{H(\curl)}}
{\|\EE_h^{\rm full}\|_{H(\curl)}}.
\end{equation}
When an analytical solution is available, the relative discretization error is
\begin{equation}\label{eq:eexnum}
e_{\rm ex}:=
\frac{\|\EE_h-\EE^{\rm ex}\|_{H(\curl)}}
{\|\EE^{\rm ex}\|_{H(\curl)}}.
\end{equation}
The numerical experiments in Examples~2--4 were implemented in MATLAB.

\subsection{Example 1: analytical verification of the symmetry-inheritance theorem on a PEC cube}
Consider $\Omega=(-1,1)^3$, centered at the origin, with a homogeneous isotropic lossless medium
\[
\mum=\I,\qquad \eps=\I,\qquad \sig=\mathbf0,
\qquad \omega=1.
\]
The isometry group of the cube is the full octahedral group $O_h$. Since the constitutive tensors are invariant under every orthogonal transformation, their symmetry group restricted to the geometry is also $O_h$, hence $G_{\rm med}=O_h$. Let
\[
a=\frac{\pi}{2},\qquad \kappa=3a^2-1,
\]
and prescribe
\begin{equation}\label{eq:ex1source}
\JJ(\xx)=j\kappa
\begin{pmatrix}
\sin(ax)\cos(ay)\cos(az)\\[1mm]
\cos(ax)\sin(ay)\cos(az)\\[1mm]
-2\cos(ax)\cos(ay)\sin(az)
\end{pmatrix}.
\end{equation}
To determine $G_J$, define
\[
\psi(\xx)=\cos(ax)\cos(ay)\cos(az),
\qquad
D=\operatorname{diag}(1,1,-2),
\]
so that
\[
\boldsymbol\Phi:=\frac{\JJ}{j\kappa}=-\frac1aD\nabla\psi.
\]
For every $g\in O_h$, the scalar field $\psi$ satisfies
$\psi(g^{-1}\xx)=\psi(\xx)$ because the elements of $O_h$ act as signed permutations of the coordinate axes. Consequently,
\[
g\nabla\psi(g^{-1}\xx)=\nabla\psi(\xx),
\]
and therefore
\[
(\TT_g^p\boldsymbol\Phi)(\xx)
=-\frac1a(gDg^T)\nabla\psi(\xx).
\]
It follows that
\[
\TT_g^p\JJ=\JJ
\quad\Longleftrightarrow\quad
gDg^T=D.
\]
The matrix $D$ has a two-dimensional eigenspace associated with the eigenvalue $1$ and the one-dimensional eigenspace $\operatorname{span}\{\mathbf e_z\}$ associated with $-2$. Hence its stabilizer in $O_h$ consists exactly of the $16$ orthogonal symmetries preserving the $z$ axis as an unoriented axis, namely $D_{4h}$. Thus
\begin{equation}\label{eq:ex1Jgroup}
G_J=D_{4h},
\qquad
G_{\rm data}=D_{4h}.
\end{equation}

The chosen frequency is nonresonant. For this PEC cube, the first positive Maxwell eigenvalue is $\pi^2/2>1=\omega^2$, so the driven problem is uniquely solvable. Theorem~\ref{thm:inherit} therefore predicts
\[
D_{4h}\subseteq G_E
\]
before the electric field is computed.

The solution can be obtained in closed form:
\begin{equation}\label{eq:ex1field}
\EE^{\rm ex}(\xx)=
\begin{pmatrix}
\sin(ax)\cos(ay)\cos(az)\\[1mm]
\cos(ax)\sin(ay)\cos(az)\\[1mm]
-2\cos(ax)\cos(ay)\sin(az)
\end{pmatrix}.
\end{equation}
On each pair of opposite faces, the two tangential components vanish because $a=\pi/2$, and therefore
\[
\nn\times\EE^{\rm ex}=\mathbf0
\qquad\text{on }\partial\Omega.
\]
Moreover,
\[
\diver\EE^{\rm ex}=0,
\qquad
\curl\curl\EE^{\rm ex}=3a^2\EE^{\rm ex},
\]
so that
\[
\curl\curl\EE^{\rm ex}-\EE^{\rm ex}
=(3a^2-1)\EE^{\rm ex}
=-j\JJ.
\]
Thus \eqref{eq:ex1field} is the unique solution. The same stabilizer calculation gives $G_E=D_{4h}$, and hence
\begin{equation}\label{eq:ex1theoremverify}
G_{\rm data}=D_{4h}=G_E.
\end{equation}

The induced mirror condition can also be verified analytically. For the reflection $R_x=\operatorname{diag}(-1,1,1)$ about $x=0$,
\[
\TT_{R_x}^p\EE^{\rm ex}=\EE^{\rm ex}.
\]
Direct differentiation gives
\[
\curl\EE^{\rm ex}
=
3a
\begin{pmatrix}
\cos(ax)\sin(ay)\sin(az)\\[1mm]
-\sin(ax)\cos(ay)\sin(az)\\[1mm]
0
\end{pmatrix}.
\]
At $x=0$, this vector is normal to the mirror plane. Since $\mum=\I$,
\[
\nn_S\times\mum^{-1}\curl\EE^{\rm ex}=\mathbf0
\qquad\text{on }x=0,
\]
which is precisely the PMC condition predicted by electric-even symmetry. No finite-element computation is involved in this example.

\subsection{Example 2: direct verification of the exact PMC and PEC symmetry boundaries}\label{sec:ex2}
The second example isolates the two parity subproblems so that the artificial boundary conditions derived in Section~\ref{sec:mirror} can be tested independently. Let $\Omega=(-1,1)^3$, $\omega=1$, and let the mirror plane be $z=0$, with
\[
R=\operatorname{diag}(1,1,-1).
\]
Choose the homogeneous anisotropic conductive medium
\begin{equation}\label{eq:constmat}
\begin{aligned}
\mum&=\operatorname{diag}(1,2,3),\\
\eps&=\operatorname{diag}(2,3,4),\\
\sig&=\operatorname{diag}(0.2,0.3,0.4).
\end{aligned}
\end{equation}
With $a=\pi/2$, define
\begin{subequations}\label{eq:3dfields}
\begin{align}
\mathbf F_+(\xx)&=
\begin{pmatrix}\cos(ay)\cos(az)\\0\\0\end{pmatrix},\\
\mathbf F_-(\xx)&=
\begin{pmatrix}0\\0\\\cos(ax)\cos(ay)\end{pmatrix}.
\end{align}
\end{subequations}
They satisfy
\[
\TT_R^p\mathbf F_+=\mathbf F_+,
\qquad
\TT_R^p\mathbf F_-=-\mathbf F_-,
\]
and both satisfy the PEC condition on the physical boundary. In addition,
\begin{align*}
\nn_S\times\mum^{-1}\curl\mathbf F_+&=\mathbf0
&&\text{on }z=0,\\
\nn_S\times\mathbf F_-&=\mathbf0
&&\text{on }z=0.
\end{align*}
Thus $\mathbf F_+$ is an exact electric-even/PMC solution and $\mathbf F_-$ is an exact electric-odd/PEC solution.

Direct differentiation yields
\begin{align*}
\curl(\mum^{-1}\curl\mathbf F_+)
&=\frac{5\pi^2}{24}\mathbf F_+,\\
\curl(\mum^{-1}\curl\mathbf F_-)
&=\frac{3\pi^2}{8}\mathbf F_-.
\end{align*}
The corresponding impressed currents are therefore prescribed as
\begin{align}
\JJ_+
&=j\left(\frac{5\pi^2}{24}-2+0.2j\right)\mathbf F_+,
\label{eq:ex2Jplus}\\
\JJ_-
&=j\left(\frac{3\pi^2}{8}-4+0.4j\right)\mathbf F_-.
\label{eq:ex2Jminus}
\end{align}
The two full-domain problems are solved separately and compared with their corresponding half-domain formulations. The relative errors $e_+$ and $e_-$ are defined by \eqref{eq:eexnum}, while $\delta_+$ and $\delta_-$ compare the full-domain discrete fields with the parity extensions of the corresponding half-domain fields. This experiment therefore tests the continuous PMC/PEC derivation and its edge-element realization without mixing the two parity classes.

\input{results/table_example2.tex}

Table~\ref{tab:ex2} shows the expected approximately first-order convergence of the lowest-order N\'ed\'elec approximation for both parity classes. At the same time, the full-domain discrete solutions and the parity extensions of the corresponding half-domain solutions agree to about machine precision. This independently confirms the PMC condition for the electric-even field and the PEC condition for the electric-odd field at the discrete level.

\subsection{Example 3: spatially varying anisotropic conductive tensors}\label{sec:ex3}
The third example tests the tensor transformation law rather than entrywise parity. Keep $\Omega=(-1,1)^3$, $R=\operatorname{diag}(1,1,-1)$, and $\omega=1$, but choose
\begin{equation}\label{eq:ex3mu}
\mum=
\begin{pmatrix}
1.40&0.18&0\\
0.18&1.80&0\\
0&0&1.25
\end{pmatrix},
\end{equation}
and the spatially varying tensors
\begin{subequations}\label{eq:ex3eps}
\begin{align}
\varepsilon_{11}&=2.4+0.20x^2,
&\varepsilon_{22}&=3.1+0.18y^2,\\
\varepsilon_{33}&=4.2+0.15z^2,
&\varepsilon_{12}&=0.10xy,\\
\varepsilon_{13}&=0.08z(1+0.15x^2),
&\varepsilon_{23}&=-0.07z(1+0.10y^2),
\end{align}
\end{subequations}
with $\varepsilon_{ji}=\varepsilon_{ij}$, and
\begin{subequations}\label{eq:ex3sig}
\begin{align}
\sigma_{11}&=0.18+0.02y^2,
&\sigma_{22}&=0.22+0.02x^2,\\
\sigma_{33}&=0.26+0.01z^2,
&\sigma_{12}&=0.012xy,\\
\sigma_{13}&=-0.020z,
&\sigma_{23}&=0.016z,
\end{align}
\end{subequations}
with $\sigma_{ji}=\sigma_{ij}$.
The matrices are uniformly positive definite on the cube and satisfy
\begin{equation}\label{eq:ex3mirror}
\begin{aligned}
R\mum R^T&=\mum,\\
R\eps(x,y,-z)R^T&=\eps(x,y,z),\\
R\sig(x,y,-z)R^T&=\sig(x,y,z).
\end{aligned}
\end{equation}
In particular, the $xz$ and $yz$ entries of \eqref{eq:ex3eps} and \eqref{eq:ex3sig} are odd in $z$. Hence material mirror symmetry is not an entrywise evenness condition; it is exactly the tensorial invariance expressed by \eqref{eq:ex3mirror}.

Use the mixed-parity analytical field
\begin{equation}\label{eq:ex3field}
\EE^{\rm ex}=0.7\mathbf F_++0.4\mathbf F_-,
\end{equation}
where $\mathbf F_\pm$ are given by \eqref{eq:3dfields}. The impressed current is manufactured from the strong equation,
\begin{equation}\label{eq:ex3source}
\JJ^{\rm ex}
=\frac{j}{\omega}
\left[
\curl(\mum^{-1}\curl\EE^{\rm ex})
-\omega^2\eps\EE^{\rm ex}
+j\omega\sig\EE^{\rm ex}
\right].
\end{equation}
The total current in \eqref{eq:ex3source} is then projected by \eqref{eq:proj}; no parity information from the known analytical field is inserted into the reduced solves. Consequently, this example checks simultaneously the tensor symmetry, the source projectors, the two artificial boundary conditions, and the reconstruction formula \eqref{eq:reconstruct}.

The full-domain and reconstructed solutions in Table~\ref{tab:ex3} exhibit essentially identical discretization errors, while the reconstruction defect remains at the level of floating-point roundoff. The observed convergence rates computed from successive refinements decrease from approximately $1.04$ to $1.01$, in agreement with the expected first-order behavior. These results also confirm that mirror invariance of a constitutive tensor is a transformed-tensor condition rather than entrywise evenness: the tangential--normal couplings are odd across the mirror plane while the complete tensors satisfy \eqref{eq:ex3mirror}.

\input{results/table_example3.tex}

\Needspace{10\baselineskip}
\subsection{Example 4: an anisotropically loaded shielded microwave cavity with an asymmetric probe}\label{sec:ex4}
The final example is intended to demonstrate both the engineering meaning and the computational consequence of the proposed reduction. Consider the normalized PEC cavity
\[
\Omega=(-1.2,1.2)\times(-0.8,0.8)\times(-0.6,0.6),
\qquad \omega=3,
\]
with mirror plane $z=0$. The background medium is
\[
\eps_b=\operatorname{diag}(1.00,1.05,0.95),
\qquad
\mum_b=\I,
\qquad
\sig_b=0.025\I.
\]
Two anisotropic lossy loading blocks are placed symmetrically on the two sides of the mirror plane. The upper block occupies
\[
B_+=\{(x,y,z):|x|\le0.60,\ |y|\le0.40,\ 0.20\le z\le0.40\},
\]
and $B_-=RB_+$. In $B_+$, set
\begin{subequations}\label{eq:ex4materials}
\begin{align}
\eps_U&=
\begin{pmatrix}
4.8&0.35&0.55\\
0.35&3.6&-0.28\\
0.55&-0.28&2.9
\end{pmatrix},\\
\mum_U&=
\begin{pmatrix}
1.18&0.06&0.08\\
0.06&1.10&-0.05\\
0.08&-0.05&1.25
\end{pmatrix},\\
\sig_U&=
\begin{pmatrix}
0.12&0.010&0.018\\
0.010&0.10&-0.014\\
0.018&-0.014&0.09
\end{pmatrix}.
\end{align}
\end{subequations}
The tensors in $B_-$ are defined by
\begin{equation}\label{eq:ex4lower}
\eps_L=R\eps_U R^T,
\qquad
\mum_L=R\mum_U R^T,
\qquad
\sig_L=R\sig_U R^T.
\end{equation}
Hence the entire piecewise anisotropic medium is mirror symmetric even though the tangential--normal couplings change sign between the two loading blocks. This configuration models a shielded microwave cavity with a pair of mirror-related anisotropic lossy inserts.

The cavity is excited by a localized off-center impressed current serving as a model of a small electric probe,
\begin{equation}\label{eq:ex4source}
\JJ(x,y,z)
=\exp[-\Phi(x,y,z)]
\begin{pmatrix}
1\\0.40j\\-0.25+0.10j
\end{pmatrix},
\end{equation}
where
\[
\Phi=
\left(\frac{x-0.46}{0.16}\right)^2
+\left(\frac{y+0.18}{0.13}\right)^2
+\left(\frac{z-0.47}{0.10}\right)^2.
\]
The excitation is deliberately placed away from the mirror plane and satisfies
\begin{equation}\label{eq:ex4noparity}
\TT_R^p\JJ\neq\JJ,
\qquad
\TT_R^p\JJ\neq-\JJ.
\end{equation}
Therefore both projected currents $\JJ_+$ and $\JJ_-$ are nonzero and both reduced problems are required. No analytical field is used in this example. The direct full-domain N\'ed\'elec solution serves only as a numerical reference for testing the exact decomposition and reconstruction.

\input{results/table_example4.tex}

\begin{figure}[!t]
\centering
\includegraphics[width=0.92\textwidth]{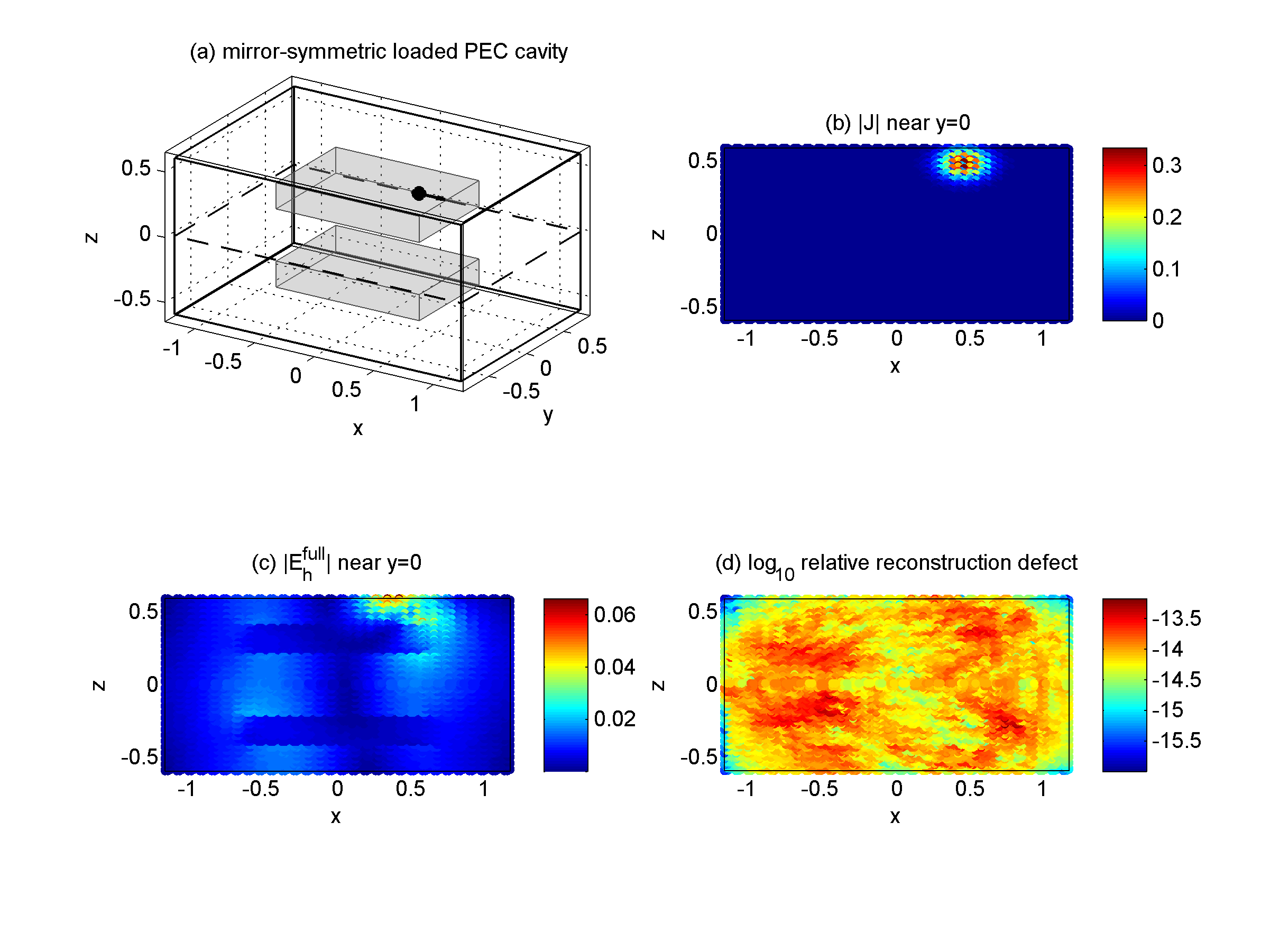}
\caption{Example 4: engineering-motivated cavity test. The two anisotropic loading blocks are mirror related, whereas the localized impressed-current model of the probe has no definite mirror parity. The field reconstructed from the PMC-even and PEC-odd half-domain solutions is compared with the direct full-domain solution.}
\label{fig:ex4overview}
\end{figure}

In Table~\ref{tab:ex4}, $\mathrm{nnz}_{LU}$ denotes the sum of the numbers of nonzero
entries in the sparse factors $\mathbf L$ and $\mathbf U$.
The quantity $\mathrm{nnz}_{LU}^{\rm full}$ refers to the full-domain
system, whereas
$\mathrm{nnz}_{LU}^{\rm halves}
=\mathrm{nnz}_{LU}^{+}+\mathrm{nnz}_{LU}^{-}$
is the combined count for the even and odd reduced systems. For the two reduced systems, the sequential time is $t_{\rm seq}=t_++t_-$, whereas $t_{\rm par}=\max(t_+,t_-)$ is reported as an ideal two-process parallel-time estimate. The factor-storage quantity $M_{\rm half}=\max(M_+,M_-)$ measures the peak storage of a single reduced factorization. These quantities distinguish the generic asymmetric-source case from the definite-parity case: although $N_++N_-=N_{\rm full}$ on the paired meshes, the method replaces one large three-dimensional factorization by two independent smaller factorizations and therefore can reduce fill-in and per-solve memory substantially.

Table~\ref{tab:ex4} demonstrates that the decomposition remains algebraically exact to near machine precision even though both parity problems are required. On the finest mesh, the direct system contains $119780$ active edge unknowns and the reconstruction defect is $1.266\times10^{-13}$. The peak sparse-LU factor storage is reduced from $4622.3$~MB for the full-domain system to $1456.5$~MB for the larger of the two reduced factorizations, a reduction of about $68.5\%$. The measured sequential time of the two reduced solves is $21.51$~s compared with $46.90$~s for the full-domain solve, while $11.37$~s is the ideal two-process estimate. 

Figure~\ref{fig:ex4overview}(a) shows the PEC cavity and the two
mirror-related anisotropic loading blocks.
Figures~\ref{fig:ex4overview}(b)--(d) show distributions on the
$x$--$z$ section at $y=0$. This observation section intersects the
mirror-symmetry plane $z=0$ along the line $z=0$.
In Figure~\ref{fig:ex4overview}(b), the current magnitude is
concentrated in the $z>0$ portion of the section, with no corresponding
concentration at the reflected locations in $z<0$.
Figure~\ref{fig:ex4overview}(c) shows that the full-domain
electric-field magnitude is also asymmetric with respect to
$z\mapsto -z$ on this section, with its largest values near the
localized excitation. These distributions illustrate that the
mirror symmetry of the geometry and medium alone does not ensure
mirror symmetry of the driven field. Figure~\ref{fig:ex4overview}(d) shows the base-10 logarithm of the
relative reconstruction defect on the same observation section.
The small values indicate close agreement, at the sampled locations,
between the direct full-domain discrete field and the field
reconstructed from the even/PMC and odd/PEC half-domain solutions.
This comparison assesses the consistency of the discrete
decomposition and reconstruction.

Examples~2--4 therefore test complementary aspects of the theory rather than repeating the same cube calculation. Example~2 directly isolates the two exact internal boundary conditions. Example~3 tests the tensor action for a spatially varying anisotropic conductive medium and verifies the full source-projection/reconstruction chain against a known analytical field. Example~4 removes the availability of a closed-form solution and uses an asymmetric probe in a mirror-symmetric, anisotropically loaded PEC cavity to demonstrate both exact discrete reconstruction and the computational consequences of splitting one full-domain system into two independent half-domain systems.

\section{Conclusion}\label{sec:conclusion}
A symmetry-inheritance and exact mirror-reduction framework has been developed directly for the second-order electric-field vector wave equation with an entirely PEC physical boundary. The derivation distinguishes polar-vector and axial-vector actions and treats the constitutive parameters as second-order tensor fields, so the reflection sign of curl and the axial character of the PEC trace are retained explicitly. Applying these actions to the original strong problem yields the transformed boundary-value problem, and unique solvability then implies that every common symmetry of the geometry, medium, boundary condition, and impressed current is inherited by the electric field.

The central result addresses source-induced mirror-symmetry breaking. For any impressed current in a mirror-symmetric background problem, reflection projectors first decompose the source on the full domain into electric-even and electric-odd components. The corresponding full-domain fields inherit these parities by uniqueness and recover the original solution by superposition. Only then is the domain reduced. The electric-odd component produces an exact PEC internal boundary, whereas the electric-even component produces an exact PMC internal boundary. The resulting half-domain problems are therefore consequences of the continuous Maxwell system rather than boundary conditions imposed heuristically from geometry alone.

Future work will extend the present framework to electromagnetic scattering problems in unbounded domains. This extension will require incorporating the outgoing radiation condition into the symmetry analysis and establishing symmetry inheritance under suitable uniqueness assumptions. Particular attention will be given to how the incident-field symmetry determines the mirror parity of the scattered field and the corresponding exact conditions on the symmetry plane. For incident fields without definite mirror parity, source decomposition and superposition will be investigated as a basis for half-domain scattering computations. Numerical implementations will also examine whether domain truncation and absorbing treatments preserve the required reflection symmetry.

\section*{Acknowledgments}
This work is supported by the National Natural Science Foundation of China under Grant No.~61901131.

\FloatBarrier
\bibliographystyle{tfcse}
\bibliography{references}

\appendix
\section{Transformation Laws for the Electric-Field Vector Wave Operator and PEC Trace}\label{app:actions}

This appendix derives the transformation formulas used in the main text. Let $g\in\OO(3)$, set $\mathbf y=g^{-1}\xx$, and define the polar-vector, axial-vector, and second-order tensor actions by \eqref{eq:mainactions}. Because $\det(g_1g_2)=\det g_1\det g_2$ and $(g_1g_2)^{-1}=g_2^{-1}g_1^{-1}$, direct substitution gives
\[
\TT_{g_1}^p\TT_{g_2}^p=\TT_{g_1g_2}^p,
\qquad
\TT_{g_1}^a\TT_{g_2}^a=\TT_{g_1g_2}^a,
\qquad
\TT_{g_1}^2\TT_{g_2}^2=\TT_{g_1g_2}^2,
\]
so all three are left actions.

The tensor action is compatible with both vector types:
\begin{align*}
(\TT_g^2\mathbf A)(\TT_g^p\mathbf v)&=\TT_g^p(\mathbf A\mathbf v),\\
(\TT_g^2\mathbf A)(\TT_g^a\mathbf w)&=\TT_g^a(\mathbf A\mathbf w),
\end{align*}
and
\[
(\TT_g^2\mathbf A)^{-1}=\TT_g^2(\mathbf A^{-1}).
\]
For two polar vectors the orthogonal cross-product identity gives
\begin{equation}\label{eq:crossaction}
(\TT_g^p\mathbf a)\times(\TT_g^p\mathbf b)
=\TT_g^a(\mathbf a\times\mathbf b).
\end{equation}
Thus the cross product of two polar vectors transforms as an axial vector.

To obtain the curl law, write $\mathbf y=g^{-1}\xx=g^T\xx$. Using the Levi--Civita symbol,
\begin{align*}
[\curl(\TT_g^p\mathbf v)]_i
&=\epsilon_{ijk}\frac{\partial}{\partial x_j}
\bigl(g_{k\ell}v_\ell(\mathbf y)\bigr)\\
&=\epsilon_{ijk}g_{k\ell}g_{jm}
\frac{\partial v_\ell}{\partial y_m}(\mathbf y)\\
&=(\det g)g_{ir}\epsilon_{rm\ell}
\frac{\partial v_\ell}{\partial y_m}(\mathbf y).
\end{align*}
Hence
\begin{equation}\label{eq:curlmain}
\curl(\TT_g^p\mathbf v)=\TT_g^a(\curl\mathbf v).
\end{equation}
Since $\TT_g^a\mathbf w=(\det g)\TT_g^p\mathbf w$, applying \eqref{eq:curlmain} once more gives
\begin{equation}\label{eq:curlaxial}
\curl(\TT_g^a\mathbf w)=\TT_g^p(\curl\mathbf w).
\end{equation}
Equations~\eqref{eq:curlmain}--\eqref{eq:curlaxial} show explicitly that curl interchanges polar and axial transformation types.

Combining tensor compatibility, tensor inversion, and the two curl identities gives
\begin{align}
\curl\!\left[(\TT_g^2\mathbf A)^{-1}
\curl(\TT_g^p\mathbf v)\right]
=\TT_g^p\!\left[\curl(\mathbf A^{-1}\curl\mathbf v)\right].
\label{eq:tensorcurlappendix}
\end{align}
Indeed, the first curl is axial, multiplication by the transformed ordinary tensor preserves the axial transformation law, and the second curl returns a polar vector. Likewise,
\[
(\TT_g^2\eps)(\TT_g^p\mathbf v)=\TT_g^p(\eps\mathbf v),
\qquad
(\TT_g^2\sig)(\TT_g^p\mathbf v)=\TT_g^p(\sig\mathbf v).
\]
Introduce the original and transformed medium data
\[
\Med=(\eps,\mum,\sig),
\qquad
\Med_g=(\TT_g^2\eps,\TT_g^2\mum,\TT_g^2\sig).
\]
The preceding identities give the operator covariance relation
\begin{equation}\label{eq:operatorcov}
\cL_{\Med_g,\omega}(\TT_g^p\mathbf v)
=\TT_g^p\!\left(\cL_{\Med,\omega}\mathbf v\right).
\end{equation}
To see explicitly how the transformed differential equation is obtained, apply
$\TT_g^p$ to the original strong equation. One obtains
\begin{align*}
\TT_g^p\!\left[
\curl(\mum^{-1}\curl\EE)
-\omega^2\eps\EE
+j\omega\sig\EE
\right]=-j\omega\TT_g^p\JJ.
\end{align*}
Using \eqref{eq:tensorcurlappendix} and tensor--vector compatibility term by term gives
\[
\curl(\mum_g^{-1}\curl\EE_g)
-\omega^2\eps_g\EE_g
+j\omega\sig_g\EE_g
=-j\omega\JJ_g,
\]
which is precisely \eqref{eq:transformed-a}. Thus the transformed equation is obtained by acting on the original PDE, rather than by postulating a transformed operator.

It remains to transform the PEC trace. Let $\nn_\Omega(\mathbf y)$ be the outward unit normal to $\partial\Omega$. Orthogonality maps tangent spaces onto tangent spaces and preserves the distinction between the interior and exterior of the mapped domain, so the outward normal to $\partial(g\Omega)$ satisfies
\begin{equation}\label{eq:normalmain}
\nn_{g\Omega}(\xx)=g\nn_\Omega(g^{-1}\xx)
=(\TT_g^p\nn_\Omega)(\xx).
\end{equation}
Using \eqref{eq:crossaction},
\begin{align*}
\nn_{g\Omega}\times(\TT_g^p\mathbf v)
=(\TT_g^p\nn_\Omega)\times(\TT_g^p\mathbf v)
=\TT_g^a(\nn_\Omega\times\mathbf v).
\end{align*}
When $g\in G_\Omega$, one has $g\Omega=\Omega$ and therefore $\TT_g^p\nn=\nn$ as an outward-normal \emph{field} on $\partial\Omega$. Consequently,
\begin{equation}\label{eq:pecappendix}
\Bt(\TT_g^p\mathbf v)
=\TT_g^a(\Bt\mathbf v),
\qquad g\in G_\Omega.
\end{equation}
This is exactly \eqref{eq:pecmain}. Applying the axial-vector action directly to the original PEC equation gives
\[
\TT_g^a(\Bt\EE)=\TT_g^a\mathbf0=\mathbf0.
\]
By \eqref{eq:pecappendix}, this is equivalent to
\[
\Bt(\TT_g^p\EE)=\mathbf0
\qquad\text{on }\partial\Omega,
\]
which is the transformed boundary condition \eqref{eq:transformed-b}. Notice that the statement $\TT_g^p\nn=\nn$ refers to the complete outward-normal field: pointwise it means $g\nn(\mathbf y)=\nn(g\mathbf y)$, not that $g\nn(\mathbf y)=\nn(\mathbf y)$ at the same boundary point.


\end{document}

%% file: results/table_example2.tex
\begin{table}[t]
\caption{Example 2: convergence and reconstruction accuracy of the electric-even/PMC and electric-odd/PEC reduced problems.}\label{tab:ex2}
\begin{tabular}{rrrrrrrrr}
\toprule
$n$ & $N_{\rm full}$ & $N_+$ & $N_-$ & $e_+$ & rate & $e_-$ & rate & $\max(\delta_+,\delta_-)$\\
\midrule
6 & 1206 & 651 & 555 & $2.683\mathrm{e}-01$ & -- & $2.103\mathrm{e}-01$ & -- & $4.454\mathrm{e}-15$\\
8 & 3032 & 1604 & 1428 & $1.838\mathrm{e}-01$ & 1.31 & $1.575\mathrm{e}-01$ & 1.01 & $1.524\mathrm{e}-14$\\
12 & 10836 & 5622 & 5214 & $1.132\mathrm{e}-01$ & 1.19 & $1.048\mathrm{e}-01$ & 1.01 & $2.014\mathrm{e}-14$\\
16 & 26416 & 13576 & 12840 & $8.230\mathrm{e}-02$ & 1.11 & $7.849\mathrm{e}-02$ & 1.00 & $9.828\mathrm{e}-15$\\
24 & 91656 & 46668 & 44988 & $5.356\mathrm{e}-02$ & 1.06 & $5.227\mathrm{e}-02$ & 1.00 & $1.317\mathrm{e}-14$\\
\bottomrule
\end{tabular}
\end{table}

%% file: results/table_example3.tex
\begin{table}[t]
\centering
\caption{Example 3: Spatially varying anisotropic conductive tensors.}
\label{tab:ex3}
\setlength{\tabcolsep}{3pt}
\renewcommand{\arraystretch}{1.05}
\begin{tabular}{rrrrrrr}
\toprule
$n$
& $N_{\rm full}$
& $N_+$
& $N_-$
& $e_{\rm full}$
& $e_{\rm rec}$
& $\delta_{\rm rec}$\\
\midrule
6
& 1206
& 651
& 555
& $2.141\mathrm{e}{-1}$
& $2.141\mathrm{e}{-1}$
& $2.591\mathrm{e}{-15}$\\
8
& 3032
& 1604
& 1428
& $1.586\mathrm{e}{-1}$
& $1.586\mathrm{e}{-1}$
& $3.004\mathrm{e}{-15}$\\
12
& 10836
& 5622
& 5214
& $1.048\mathrm{e}{-1}$
& $1.048\mathrm{e}{-1}$
& $7.163\mathrm{e}{-15}$\\
16
& 26416
& 13576
& 12840
& $7.835\mathrm{e}{-2}$
& $7.835\mathrm{e}{-2}$
& $7.662\mathrm{e}{-15}$\\
24
& 91656
& 46668
& 44988
& $5.210\mathrm{e}{-2}$
& $5.210\mathrm{e}{-2}$
& $1.557\mathrm{e}{-14}$\\
\bottomrule
\end{tabular}
\end{table}

%% file: results/table_example4.tex
\begin{table}[t]
\centering
\caption{Example 4: computational performance for the anisotropically loaded PEC cavity. }
\label{tab:ex4}
\setlength{\tabcolsep}{4pt}
\textit{(a) System sizes and sparse-LU factor storage}\par\smallskip
\begin{tabular}{crrrcc}
\toprule
Mesh & $N_{\rm full}$ & $N_+$ & $N_-$ & $\mathrm{nnz}_{LU}^{\rm full}/\mathrm{nnz}_{LU}^{\rm halves}$ & $M_{\rm full}/M_{\rm half}$ (MB)\\
\midrule
12$\times$8$\times$12 & 7088 & 3678 & 3410 & 2.428e+06/1.448e+06 & 56.0/17.5\\
18$\times$12$\times$18 & 24996 & 12807 & 12189 & 1.577e+07/1.040e+07 & 362.5/127.7\\
24$\times$16$\times$24 & 60544 & 30828 & 29716 & 6.392e+07/4.381e+07 & 1466.8/538.2\\
30$\times$20$\times$30 & 119780 & 60765 & 59015 & 2.016e+08/1.216e+08 & 4622.3/1456.5\\
\bottomrule
\end{tabular}
\par\medskip
\textit{(b) Solve times and reconstruction defects}\par\smallskip
\begin{tabular}{crrrc}
\toprule
Mesh & $t_{\rm full}$ (s) & $t_{\rm seq}$ (s) & $t_{\rm par}$ (s) & $\delta_{\rm rec}$\\
\midrule
12$\times$8$\times$12 & 0.25 & 0.16 & 0.08 & 6.663e-14\\
18$\times$12$\times$18 & 2.23 & 1.12 & 0.63 & 1.123e-13\\
24$\times$16$\times$24 & 11.72 & 6.18 & 3.55 & 1.122e-13\\
30$\times$20$\times$30 & 46.90 & 21.51 & 11.37 & 1.266e-13\\
\bottomrule
\end{tabular}
\end{table}